\documentclass[11pt,a4paper,oneside]{article}

\usepackage{a4,amsthm,amsfonts}
\usepackage[intlimits]{amsmath}

\author{Stefan Boller\\
\\
Graduiertenkolleg Quantenfeldtheorie\\
Universit\"at Leipzig, Augustusplatz 10/11\\
04109 Leipzig, Germany\\
e-mail: Stefan.Boller@itp.uni-leipzig.de}

\title{Characterization of Cyclic and Separating Vectors and Application to an
  Inverse Problem in Modular Theory\\
  II. Semifinite Factors}

\newtheorem{thm}{Theorem}[section]
\newtheorem{cor}[thm]{Corollary}
\newtheorem{lem}[thm]{Lemma}
\newtheorem{prop}[thm]{Proposition}

\theoremstyle{definition}
\newtheorem{defn}{Definition}[section]

\theoremstyle{remark}
\newtheorem{rem}{Remark}[section]

\newtheorem{bsp}{Example}[section]

\numberwithin{equation}{section}

\newcommand{\thmref}[1]{Theorem~\ref{#1}}
\newcommand{\secref}[1]{\S\ref{#1}}
\newcommand{\lemref}[1]{Lemma~\ref{#1}}
\newcommand{\corref}[1]{Corollary~\ref{#1}}
\newcommand{\propref}[1]{Proposition~\ref{#1}}
\newcommand{\defnref}[1]{Definition~\ref{#1}}
\newcommand{\remref}[1]{Remark~\ref{#1}}

\newcommand{\bspref}[1]{Example~\ref{#1}}

\DeclareMathOperator{\RE}{Re}

\DeclareMathOperator{\tr}{tr}
\DeclareMathOperator{\lin}{lin}

\newcommand{\interval}[1]{\mathinner{#1}}
\newcommand{\offen}[2]{\interval{(#1,#2)}}
\newcommand{\abge}[2]{\interval{\lbrack#1,#2\rbrack}}

\newcommand{\eval}[2][\right]{\relax
  \ifx#1\right\relax \left.\fi#2#1\rvert}

\let\norm=\enVert

\newcommand{\spitz}[1]{\bigl\langle#1\bigr\rangle}
\newcommand{\SProd}[2]{\spitz{#1\rvert#2}}

\newcommand{\field}[1]{\mathbb{#1}}
\newcommand{\R}{\field{R}}
\newcommand{\C}{\field{C}}

\newcommand{\N}{\field{N}}

\newcommand{\raum}[1]{\mathcal{#1}}
\renewcommand{\H}{\raum{H}}

\newcommand{\alg}[1]{\mathcal{#1}}
\newcommand{\M}{\alg{M}}

\newcommand{\Op}[1]{\mathrm{#1}}

\newcommand{\Dom}[1]{\mathcal{D}(\Op{#1})}

\newcommand{\Lim}[1]{\lim_{#1\rightarrow\infty}}

\begin{document}

  \maketitle

  \begin{abstract} 
    This paper generalizes the results obtained in an earlier paper
    (\cite{BolI}) for finite
    factors to infinite but still semifinite factors.
    First we give a characterization of cyclic and
    separating vectors for infinite semifinite factors in terms of operators
    associated with this vector and being
    affiliated with the factor. Further we show how this operator generates
    the 
    modular objects of the given cyclic and separating vector generalizing an
    idea of Kadison and Ringrose. With the help
    of these results we can show that the second simple class of solutions for
    the inverse problem constructed in \cite{BolI} never exists in infinite
    semifinite factors. Finally we give a
    classification of the solutions of the inverse problem in the case of
    modular operators having pure point
    spectrum completely analoguous to the finite case. 
  \end{abstract}

  \setcounter{page}{1}

\section{The Inverse Problem in Modular Theory}
\label{sec:1}

Let $\M_{0}$ be a von Neumann algebra on a separable Hilbert space $\H_{0}$
with a cyclic and separating vector $u_{0}$. Then modular theory shows the
existence of a modular operator $\Delta_{0}$ and a modular conjugation
$\Op{J}_{0}$ (the modular objects $(\Delta_{0},\Op{J}_{0})$) belonging to the
vector $u_{0}$. In this paper we examine the inverse problem of constructing
algebras $\M$ having the same cyclic and separating vector and modular objects
as $\M_{0}$:

\textbf{The Inverse Problem}

Let $(\Delta_{0},\Op{J}_{0})$ be the modular objects for the von Neumann
algebra $\M_{0}$ with cyclic and separating vector $u_{0}$. Characterize all
von Neumann algebras $\M$ isomorphic to $\M_{0}$ with the following
properties:
\begin{enumerate}
  \item $u_{0}$ is also cyclic and separating for $\M$,
  \item $(\Delta_{0},\Op{J}_{0})$ are the modular objects for $(\M,u_{0})$.
\end{enumerate}
Let $NF_{\M_{0}}(\Delta_{0},\Op{J}_{0},u_{0})$ denote all solutions $\M$ of
the inverse problem.

In \cite{BolI} the following theorems were shown:
\begin{thm}\label{thm1:1}
  Let ($\M_{0},\H_{0}$) be a finite von Neumann factor. Let further
  $u\in\H_{0}$. Then there is
  exactly one operator $\Op{T}_{u}\eta\M_{0}$ associated with the vector $u$,
  s.t. $u=\Op{T}_{u}u_{\tr}$ where $u_{\tr}\in\H_{0}$ is a cyclic trace
  vector. This operator has the
  following properties: 
  \begin{enumerate}
    \item $\tr(\Op{T}_{u}\Op{T}_{u}^{*})=\tr(\Op{T}_{u}^{*}\Op{T}_{u})<\infty$.
    \item $u$ is cyclic, iff $\Op{T}_{u}$ is injective.
    \item $u$ is separating, iff $\Op{T}_{u}$ has dense range.
    \item $u$ is cyclic and separating iff $\Op{T}_{u}$ is injective and has
  dense range, i.e. iff $\Op{T}_{u}$ is invertible. 
  \end{enumerate}
\end{thm}

\begin{thm}\label{thm1:2}
  Let $\Op{T}\eta\M_{0}$. Then there is a vector $u\in\H_{0}$
  s.t $\Op{T}=\Op{T}_{u}$ in the sense of \thmref{thm1:1} iff
  $\tr{\Op{T}\Op{T}^{*}}=\tr{\Op{T}^{*}\Op{T}}<\infty$.  
\end{thm}

In the second section of this paper these theorems were generalized to
infinite semifinite factors. For this purpose we first consider a special case
of such factors, a matrix algebra of finite factors, where the trace vector is
replaced by a sequence of vectors, constructed from the trace vectors of the
constituting factors.

With the help of this result we show the analogue of the following result,
also obtained in \cite{BolI} for finite factors:
\begin{thm}\label{thm1:3}
  Let $\M_0$ be a finite von Neumann factor with cyclic and separating vector
  $u_{0}\in\M_{0}$ and cyclic trace vector $u_{\tr}\in\H_0$. Let further
  $\Op{T}_{u_{0}}\eta\M_{0}$ be the invertible operator corresponding to
  $u_{0}$ and
  $\Op{T}_{u_{0}}=\Op{HV}=(\Op{T}_{u_{0}}\Op{T}_{u_{0}}^{*})^{1/2}\Op{V}$ the
  polar decomposition of $\Op{T}_{u_{0}}$. Then we can calculate the modular
  objects $(\Delta_{0},\Op{J}_{0})$ of $(\M_{0},u_{0})$ as follows:
  \begin{equation*}
    \Op{J}_0=\Op{JV}^*\Op{JVJ}=\Op{VJV}^*,
  \end{equation*}
  where $\Op{J}$ is the
    conjugation corresponding to $u_{\tr}$, and
  \begin{equation*}
    \Delta_0=\Op{J}_0\Op{H}^{-1}_0\Op{J}_0\Op{H}_0,
  \end{equation*}
  where $\Op{H}_0=\Op{H}^2=\Op{T}_{u_{0}}\Op{T}_{u_{0}}^{*}$.
\end{thm}

Then we will be in exactly the same situation as in the finite case, and can
examine the inverse problem as in \cite{BolI}. In contrast to that case the
second simple class of solutions will never exists in this case
(s. \secref{sec:4}), but the classification of the solutions in the pure point
spectrum case will be the same.

Notice that in this paper all Hilbert spaces are separable, i.e. the von
Neumann algebras are countably decomposable.

\section{Characterization of Vectors by Affiliated Operators}
\label{sec:2}

In this section we consider infinite but still semifinite factors,
i.e. $\M_{0}$ is of type $I_{\infty}$ or $II_{\infty}$. In this case we have
no trace vectors left. But nevertheless we can make a similar construction as
in the finite case by considering the infinite factor as an infinite
matrix of finite factors and using the results presented in \cite{BolI}.

As a model for such a matrix of finite factors we examine now the semifinite
factor $\alg{R}=\alg{T}\otimes L(\H_{\infty})\otimes\C$ on
$\raum{K}=\H\otimes\H_{\infty}\otimes\H_{\infty}$, where $(\alg{T},\H)$ is a
finite factor possessing a cyclic and separating vector and $\H_{\infty}$ is a
infinite dimensional separable Hilbert space which we can identify with
$l_{2}(\N)$. Now $\alg{R}$ is an (infinite) type $I$ ($II$) factor, if
$\alg{T}$ is type $I$ ($II$). Further, since 
$(\alg{T},\H)$ is finite and possesses a cyclic and separating vector, it
possesses a cyclic trace vector $u_{\tr}\in\H$ (cf. \cite[Th. 8.2.8,
Lem. 7.2.8]{KRII}). 

In the following we consider the elements of $\raum{K}$ as infinite
dimensional matrices $u=(u_{i}^{k})_{i,k\in\N}$ with entries $u_{i}^{k}\in\H$
s.t. $\sum_{i,k}\norm{u_{i}^{k}}^{2}<\infty$, where the lower index
corresponds to the second component of the tensor product and the upper to the
third, resp. Then we
can write the elements of $\alg{R}$ as matrices
$\Op{T}=(\Op{T}_{li})_{l,i\in\N}$ with entries
$\Op{T}_{li}\in\alg{T}$, where
\begin{equation*}
  \Op{T}u=(\sum_{i}\Op{T}_{li}u_{i}^{k})_{l}^{k}.
\end{equation*}
Then the commutant $\alg{R}^{'}$ of $\alg{R}$ is
$\alg{T}^{'}\otimes\C\otimes L(\H_{\infty})$, where we can write an
element in $\alg{R}^{'}$ as $\Op{T}^{'}=({\Op{T}^{'}}^{lk})_{l,k\in\N}$ with
entries
${\Op{T}^{'}}^{lk}\in\alg{T}^{'}$, where
\begin{equation*}
  \Op{T}^{'}u=(\sum_{k}{\Op{T}^{'}}^{lk}u_{i}^{k})_{i}^{l}.
\end{equation*}

For the proofs in this section the following subalgebras of $\alg{R}$ and
$\alg{R}^{'}$ are important:
\begin{equation*}
  \begin{split}
    \alg{R}_{0}&:=\{\Op{M}=(\Op{M}_{ij})_{ij}\in\alg{R}\vert
    \Op{M}_{ij}\not=0\text{ for only finitely many }i,j\in\N\}\\
    \alg{R}_{0}^{'}&:=\{\Op{M}^{'}=({\Op{M}^{'}}^{ij})^{ij}\in\alg{R}^{'}\vert
    {\Op{M}^{'}}^{ij}\not=0\text{ for only finitely many }i,j\in\N\}.
  \end{split}
\end{equation*}
Now we define the following sequence of vectors in $\raum{K}$, which is
the analogue to the trace vector:
\begin{equation*}
  v_{k}:=(\delta_{i}^{j}\delta_{ik}u_{\tr})_{i}^{j}.
\end{equation*}
Further we define 
\begin{equation*}
  D_{0}:=\lin\{\Op{M}v_{k}\vert
    \Op{M}\in\alg{R},k\in\N\}\subset\raum{K}
\end{equation*}
and
\begin{equation*}
  D_{0}^{'}:=\lin\{\Op{M}^{'}v_{k}\vert
    \Op{M}^{'}\in\alg{R}^{'},k\in\N\}\subset\raum{K}.
\end{equation*}

Now the trace $\tr$ of $\alg{R}$, which is a n.s.f. tracial weight, is
$\tr=\tr_{\alg{T}}\otimes\tr_{L(\H_{\infty})}$, where $\tr_{\alg{T}}$ is the
trace on $\alg{T}$ and $\tr_{L(\H_{\infty})}$ the standard trace on
$\H_{\infty}$. It can be written with the help of the vectors $(v_{k})$:
\begin{equation*}
  \begin{split}
    \tr(\Op{M})
    &=\sum_{k}\tr_{\alg{T}}(\Op{M}_{kk})
    =\sum_{k}\SProd{\Op{M}_{kk}u_{\tr}}{u_{\tr}}\\
    &=\sum_{k}\SProd{\Op{M}v_{k}}{v_{k}}\\
    &=\sum_{k}\int\lambda d\norm{\Op{E}_{\lambda}^{\Op{M}}v_{k}}^{2}
    \quad\forall\Op{M}=(\Op{M}_{ij})\in\alg{R}^{+},
  \end{split}
\end{equation*}
where $\Op{M}=\int\lambda d\Op{E}_{\lambda}^{\Op{M}}$ is the spetral measure
of $\Op{M}$.
As in \cite{BolI} we can continue the trace to all the positive closed
operators $\Op{A}$ affiliated  with $\alg{R}$ by 
\begin{equation}\label{eq2:2:1}
  \tr(\Op{A}):=\sum_{k}\int\lambda d\norm{\Op{E}_{\lambda}^{\Op{A}}v_{k}}^{2},
\end{equation}
where $\Op{E}_{\lambda}^{\Op{A}}$ is the spectral measure of $\Op{A}$.

Now we can associate an operator $\Op{T}_{ij}\eta\alg{T}$ with every component
$u_{i}^{j}$ of a vector $u=(u_{i}^{j})_{i}^{j}\in\H$,
s.t. $u_{\tr}\in\Dom{T_{ij}}$ and $\Op{T}_{ij}u_{\tr}=u_{i}^{j}$
(cf. \cite{BolI}). These operators give rise to a linear operator
$\tilde\Op{T}_{u}$ defined by
\begin{equation}\label{eq2:2:2}
  \begin{split}
    \tilde\Op{T}_{u}:\Dom{\tilde T_{u}}:=D_{0}^{'}\subset\raum{K}&\to\raum{K}\\
    \Op{M}^{'}v_{k}&\mapsto\tilde\Op{T}_{u}\Op{M}^{'}v_{k}
    :=(\Op{T}_{ik}{\Op{M}^{'}}^{jk}u_{\tr})_{i}^{j}.
  \end{split}
\end{equation}
Now we can prove
\begin{lem}\label{lem2:2:2a}
  Let $\tilde\Op{T}_{u}$ defined by \eqref{eq2:2:2}. Then $\tilde\Op{T}_{u}$
  is densely defined and closable. Let $\Op{T}_{u}$ be its closure. Then
  $D_{0}^{'}\subset\Dom{\Op{T}_{u}^{*}}$, $\Op{T}_{u}$ is affiliated with
  $\alg{R}$, and
  \begin{equation*}
    \sum_{k}\Op{T}_{u}v_{k}=u,
  \end{equation*}
  where the convergence is absolute.
\end{lem}

\begin{proof}
  \begin{enumerate}
  \item
  First we must show that $\tilde\Op{T}_{u}$ is well defined. Observe first
  that
  \begin{equation*}
    \begin{split}
      \sum_{i,j}\norm{\Op{T}_{ik}\delta^{jk}u_{\tr}}^{2}
      &=\sum_{i}\norm{\Op{T}_{ik}u_{\tr}}^{2}\\
      &=\sum_{i}\norm{u_{i}^{k}}^{2}<\infty,
    \end{split}
  \end{equation*}
  i.e. $\tilde\Op{T}_{u}v_{k}\in\raum{K}$ for every $k\in\N$. Let now
  $\Op{M}^{'}\in\alg{R}^{'}$ be arbitrary. Then
  $\Op{M}^{'}\tilde\Op{T}_{u}v_{k}\in\raum{K}$ and
  \begin{equation}\label{eq2:2:0}
    \begin{split}
      \infty&>\norm{\Op{M}^{'}\tilde\Op{T}_{u}v_{k}}^{2}\\
      &=\sum_{i,j}\norm{{\Op{M}^{'}}^{jk}\Op{T}_{ik}u_{\tr}}^{2}\\
      &=\sum_{i,j}\norm{\Op{T}_{ik}{\Op{M}^{'}}^{jk}u_{\tr}}^{2}
      \text{ (cf. \cite[Prop.2.1.]{BolI})}\\
      &=\norm{\tilde\Op{T}_{u}\Op{M}^{'}v_{k}}^{2}
    \end{split}
  \end{equation}
  for every $k\in\N$, hence
  $\tilde\Op{T}_{u}$ is well defined.

  \item
  Now we show that $\Dom{\tilde T_{u}}$ is dense in $\raum{K}$. First the
  elements with only finitely many entries not $0$ are dense in
  $\raum{K}$. Further every such element is a linear combination of elements
  of the type $(u_{i}^{j}\delta_{ik})_{i}^{j}$, again all but a finite number
  equal $0$. Since $u_{\tr}\in\H$ is cyclic for
  $\alg{T}^{'}$, we can approximate these elements by elements of the form
  $({\Op{M}^{'}}^{jk}\delta_{ik}u_{\tr})_{i}^{j}=:\Op{M}^{'}v_{k}$ with
  $\Op{M}^{'}=({\Op{M}^{'}}^{jk}\delta^{ik})^{ji}
   \in\alg{R}_{0}^{'}\subset\alg{R}^{'}$, hence
  $\Dom{\tilde T_{u}}=D_{0}^{'}$ is dense in $\raum{K}$. 

  \item
  In this step we want to show that $\tilde\Op{T}_{u}$ is closable. Let
  $x=\Op{M}^{'}v_{k}\in\Dom{\tilde T_{u}}$,
  $y=\Op{N}^{'}v_{j}\in\Dom{S}:=D_{0}^{'}$ ($k,j\in\N$),
  where $\Op{S}:=(\Op{T}^{*}_{li})_{i,l}$ is defined analogously to
  $\tilde\Op{T}_{u}$, hence it is a densely defined operator, too (All
  $\Op{T}_{li}^{*}$ are closed operators affiliated with $\alg{T}$ and
  $u_{\tr}\in\Dom{T_{li}^{*}}$, cf. \cite[Prop. 2.1.]{BolI}, and
  $\norm{\Op{T}_{li}^{*}u_{\tr}}^{2}=\norm{\Op{T}_{li}u_{\tr}}^{2}$).

  Now
  \begin{equation*}
    \begin{split}
      \SProd{\tilde\Op{T}_{u}x}{y}
      &=\SProd{\tilde\Op{T}_{u}\Op{M}^{'}v_{k}}{\Op{N}^{'}v_{j}}\\
      &=\sum_{i,l}\SProd{\Op{T}_{ik}{\Op{M}^{'}}^{lk}u_{\tr}}
                        {{\Op{N}^{'}}^{lj}\delta_{ij}u_{\tr}}\\
      &=\sum_{l}\SProd{{\Op{M}^{'}}^{lk}u_{\tr}}
                        {\Op{T}_{jk}^{*}{\Op{N}^{'}}^{lj}u_{\tr}}\\
      &=\sum_{i,l}\SProd{\delta_{ik}{\Op{M}^{'}}^{lk}u_{\tr}}
                        {\Op{T}_{ji}^{*}{\Op{N}^{'}}^{lj}u_{\tr}}\\
      &=\SProd{x}{\Op{S}y}.
    \end{split}
  \end{equation*}
  This shows $y\in\Dom{{\tilde T_{u}}^{*}}$, $\tilde\Op{T}_{u}^{*}y=\Op{S}y$,
  and $\Op{S}\subset(\tilde\Op{T}_{u})^{*}$, hence $\tilde\Op{T}_{u}$ is
  closable. This shows also, that $D_{0}^{'}\subset\Dom{(\tilde
  T_{u})^{*}}=\Dom{T_{u}^{*}}$. 

  \item
  To show that $\Op{T}_{u}$ is affiliated with $\alg{R}$, let
  $\Op{U}^{'}=({\Op{U}^{'}}^{ij})_{i,j\in\N}\in\alg{R}^{'}$ be a unitary. 
  Then $\Op{U}^{'}D_{0}^{'}=D_{0}^{'}$.
  Let now 
  $x=\Op{M}^{'}v_{k}\in D_{0}^{'}=\Dom{\tilde T_{u}}$. Then
  \begin{equation*}
    \begin{split}
      \Op{U}^{'}\tilde\Op{T}_{u}x
      &=\Op{U}^{'}(\Op{T}_{ik}{\Op{M}^{'}}^{jk}u_{\tr})_{i}^{j}\\
      &=(\sum_{j}{\Op{U}^{'}}^{lj}\Op{T}_{ik}
                 {\Op{M}^{'}}^{jk}u_{\tr})_{i}^{l}\\
      &=(\sum_{j}\Op{T}_{ik}{\Op{U}^{'}}^{lj}
                 {\Op{M}^{'}}^{jk}u_{\tr})_{i}^{l}
          \text{ (cf. \cite[Prop.2.1.]{BolI})}\\
      &=(\Op{T}_{ik}\sum_{j}{\Op{U}^{'}}^{lj}
                 {\Op{M}^{'}}^{jk}u_{\tr})_{i}^{l}\\
      &=\tilde\Op{T}_{u}\Op{U}^{'}x
    \end{split}
  \end{equation*}
  This shows $\Op{U}^{'}\tilde\Op{T}_{u}=\tilde\Op{T}_{u}\Op{U}^{'}$ for every
  unitary $\Op{U}^{'}\in\alg{R}^{'}$, hence, since $D_{0}^{'}$ is a core for
  $\Op{T}_{u}$,
  \begin{equation*}
    \Op{U}^{'}\Op{T}_{u}=\Op{T}_{u}\Op{U}^{'}
    \quad\forall\Op{U}^{'}\in\alg{U}(\alg{R}^{'}).
  \end{equation*}
  
  \item
    In the last step we calculate
    \begin{equation*}
      \begin{split}
        \sum_{k}\Op{T}_{u}v_{k}&=
        \sum_{k}(\Op{T}_{ik}\delta^{jk}u_{tr})_{i}^{j}\\
        &=(\Op{T}_{ij}u_{tr})_{i}^{j}=u.
      \end{split}
    \end{equation*}
  \end{enumerate}
\end{proof}

Now we can give the following definition:
\begin{defn}\label{defn2:2:1}
  For every vector $u=(u_{i}^{j})\in\raum{K}$ we denote by
  $\Op{T}_{u}=(\Op{T}_{ij})$ an
  operator affiliated with $\alg{R}$ s.t. $u_{\tr}\in\Dom{T_{ij}}$ for all
  $i,j\in\N$, $\Op{T}_{ij}u_{\tr}=u_{i}^{j}$, and $\sum_{k}\Op{T}_{u}v_{k}=u$,
  which exists according to \lemref{lem2:2:2a}.
\end{defn}

The next proposition shows some usefull properties of the operators occuring
in \defnref{defn2:2:1}  
\begin{prop}\label{prop2:2:2}
  Let $\Op{T}\eta\alg{R}$, $v_{k}\in\Dom{T}$ ($k\in\N$),
  $\sum_{k}\norm{\Op{T}v_{k}}^{2}<\infty$. Then
  \begin{enumerate}
    \item $\alg{R}v_{k}\subset\Dom{T}$,  $\alg{R}v_{k}\subset\Dom{T^{*}}$, and
      $\alg{R}v_{k}\subset\Dom{(T^{*}T)^{1/2}}$ for all $k\in\N$.
    \item $D_{0}$ is a core for $\Op{T}$, $\Op{T}^{*}$, and
      $(\Op{T^{*}}\Op{T})^{1/2}$.
    \item $\alg{R}^{'}v_{k}\subset\Dom{T}$,
      $\alg{R}^{'}v_{k}\subset\Dom{T^{*}}$, and
      $\alg{R}^{'}v_{k}\subset\Dom{(T^{*}T)^{1/2}}$ for all $k\in\N$.
    \item $D_{0}^{'}$ is a core for $\Op{T}$, $\Op{T}^{*}$, and
      $(\Op{T^{*}}\Op{T})^{1/2}$.
  \end{enumerate}
\end{prop}
\begin{proof}
  \begin{enumerate}
  \item Let $\Op{T}=\Op{V}\Op{H}$ the polar decomposition of $\Op{T}$, and
    $\Op{E}_{\lambda}$ the spectral resolution of $\Op{H}$. Then
    \begin{equation*}
      \begin{split}
        \int\lambda^{2}d\norm{\Op{E}_{\lambda}\Op{U}v_{k}}^{2}
        &\le\sum_{l}\int\lambda^{2}d\norm{\Op{E}_{\lambda}\Op{U}v_{l}}^{2}\\
        &=\sum_{l}\int\lambda^{2}d\norm{\Op{U}^{*}\Op{E}_{\lambda}v_{l}}^{2}
        \text{(s. \eqref{eq2:2:1})}\\
        &=\sum_{l}\int\lambda^{2}d\norm{\Op{E}_{\lambda}v_{l}}^{2}\\
        &=\sum_{l}\norm{\Op{H}v_{l}}^{2}=\sum_{l}\norm{\Op{T}v_{l}}^{2}<\infty
      \end{split}
    \end{equation*}
    for every unitary $\Op{U}\in\alg{R}$ and every $k\in\N$,
    i.e. $\alg{R}v_{k}\subset\Dom{H}=\Dom{T}$ for every $k\in\N$. Now
    $\Op{T}^{*}=\Op{H}\Op{V}^{*}$, and, since $\Op{V}^{*}\in\alg{R}$, also
    $\alg{R}v_{k}\subset\Dom{T^{*}}$.
  \item 1) shows that $D_{0}\subset\Dom{T}$, further $D_{0}$ is dense in
    $\raum{K}$. Now $D_{0}$ is invariant under the unitary group
    $e^{it\Op{H}}$, i.e. $D_{0}$ is a core for $\Op{H}$ and also for
    $\Op{T}$. The assertion for $\Op{T}^{*}$ follows analogous.  
  \item This follows from 1) and \cite[Prop. 2.1.]{BolI}.
  \item Now for every $\Op{M}=(\Op{M}_{ij})_{ij}\in\alg{R}$ there exists
    exactly one $\Op{M}^{'}=({\Op{M}^{'}}^{ij})\in\alg{R}^{'}$
    s.t. ${\Op{M}^{'}}^{ij}u_{\tr}=\Op{M}_{ji}u_{\tr}$
    (${\Op{M}^{'}}^{ij}:=\Op{J}\Op{M}_{ji}\Op{J}$, where $\Op{J}$ is the
    conjugation w.r.t. $u_{\tr}$). Now define
    \begin{equation*}
      \Op{M}^{'}_{(k,l)}:=\Op{E}^{'}_{k}\Op{M}^{'}\Op{E}^{'}_{l},
    \end{equation*}
    where $\Op{E}_{(k)}^{'}:=(\delta_{ik}\delta^{ij})^{ij}$. Then
    \begin{equation}\label{eq2:2:4}
      \Op{M}v_{k}=\sum_{l}\Op{M}^{'}_{(k,l)}v_{l}
    \end{equation}
    and
    \begin{equation*}
      \begin{split}
        \sum_{l}\norm{\Op{T}\Op{M}^{'}_{(k,l)}v_{l}}^{2}
        &=\sum_{l}\norm{\Op{M}^{'}_{(k,l)}\Op{T}v_{l}}^{2}\\
        &\le\sum_{l}\norm{\Op{M}^{'}}^{2}\norm{\Op{T}v_{l}}^{2}\\
        &\le\norm{\Op{M}}^{2}\sum_{l}\norm{\Op{T}v_{l}}^{2}<\infty.
      \end{split}
    \end{equation*}
    Hence $\sum_{l}\Op{T}\Op{M}^{'}_{(k,l)}v_{l}$ converges and therefore
    $\sum_{l}\Op{M}^{'}_{(k,l)}v_{l}$ converges in the graph norm
    of $\Op{T}$ to $\Op{M}v_{k}$, i.e. also $D_{0}^{'}$ is a core, since
    $D_{0}$ is it.
  \end{enumerate}
\end{proof}

\begin{lem}\label{lem2:2:3}
  Let $\Op{T}\eta\alg{R}$ be as in \propref{prop2:2:2}. Then there are
    $\Op{T}_{ij}\eta\alg{T}$ with $u_{\tr}\in\Dom{T_{ij}}$
    s.t. $\Op{T}=(\Op{T}_{ij})_{i,j\in\N}$ and $\Op{T}=\Op{T}_{u}$ with
    $u:=\sum_{k}\Op{T}v_{k}$ in the sense of \defnref{defn2:2:1}.
\end{lem}

\begin{proof}
  Set $\Op{E}_{k}:=[\alg{R}^{'}v_{k}]\in\alg{R}$ ($k\in\N$). Then matrix
  calculation shows that
  $(\Op{E}_{k})$ is a family of orthogonal, equivalent, finite projections,
  s.t. 
  \begin{equation*}
    \Op{E}_{k}\alg{R}\Op{E}_{k}=
    (\delta_{ki}\delta_{ij}\alg{T})_{i,j\in\N} 
  \end{equation*}
  and $\sum_{k}\Op{E}_{k}=\Op{Id}$. Now
  $\Op{E}_{ij}:\Op{M}^{'}v_{j}\mapsto\Op{M}^{'}v_{i}$ defines a selfadjoint system
  of matrix units $(\Op{E}_{ij})$ s.t. $\Op{E}_{kk}=\Op{E}_{k}$ for every
  $k\in\N$. Now define operators 
  \begin{equation*}
    \begin{split}
      \Op{S}_{ij}:\Dom{S_{ij}}:=D_{0}^{'}\subset\raum{K}
      &\to\raum{K}\\
      \sum_{k}\Op{M}^{'}_{k}v_{k}&\mapsto
      \sum_{k}\Op{E}_{ki}\Op{T}\Op{E}_{jk}\Op{M}^{'}_{k}v_{k}
    \end{split}
  \end{equation*}
  and
  \begin{equation*}
    \begin{split}
      \tilde\Op{S}_{ji}:\Dom{\tilde S_{ij}}:=D_{0}^{'}\subset\raum{K}
      &\to\raum{K}\\
      \sum_{k}\Op{M}^{'}_{k}v_{k}&\mapsto
      \sum_{k}\Op{E}_{kj}\Op{T}^{*}\Op{E}_{ik}\Op{M}^{'}_{k}v_{k}.
    \end{split}
  \end{equation*}
  Since $D_{0}^{'}$ is dense in $\raum{K}$ (cf. proof of
  \lemref{lem2:2:2a}) and a core both for $\Op{T}$ and for $\Op{T}^{*}$ they
  are well defined and densely defined. Let now $x\in\Dom{S_{ij}}$ and
  $y\in\Dom{\tilde S_{ji}}$. Then
  \begin{equation*}
    \begin{split}
      \SProd{\Op{S}_{ij}x}{y}
      &=\sum_{k}\SProd{\Op{E}_{ki}\Op{T}\Op{E}_{jk}x}{y}\\
      &=\sum_{k}\SProd{x}{\Op{E}_{kj}\Op{T}^{*}\Op{E}_{ik}y}\\
      &=\SProd{x}{\Op{\tilde S}_{ji}y}.
    \end{split}
  \end{equation*}
  This means that $y\in\Dom{S_{ij}^{*}}$ and $\Op{S}_{ij}^{*}y=\Op{\tilde
    S}_{ji}y$, i.e.
  $\Op{\tilde S}_{ji}\subset\Op{S}_{ij}^{*}$, hence $\Op{S}_{ij}$ is
  closable since $\Op{\tilde S}_{ji}$ is densely defined. 
  
  Let now $\tilde\Op{T}_{ij}$ be the closure of $\Op{S}_{ij}$. Then
  $D_{0}^{'}=\Dom{S_{ij}}$ is a core for $\tilde\Op{T}_{ij}$. Since
  $\Op{U}^{'}D_{0}^{'}=D_{0}^{'}$ and
  \begin{equation*}
    \begin{split}
      \Op{U}^{'}\Op{S}_{ij}(\sum_{k}\Op{M}_{k}^{'}v_{k})
      &=\sum_{k}\Op{U}^{'}\Op{E}_{ki}\Op{T}\Op{E}_{jk}\Op{M}_{k}^{'}v_{k}\\
      &=\sum_{k}\Op{E}_{ki}\Op{T}\Op{E}_{jk}\Op{U}^{'}\Op{M}_{k}^{'}v_{k}\\
      &=\Op{S}_{ij}\Op{U}^{'}(\sum_{k}\Op{M}_{k}^{'}v_{k})
    \end{split}
  \end{equation*}
  for every unitary $\Op{U}^{'}\in\alg{R}^{'}$ and every element
  $(\Op{M}^{'}_{k})v_{k}\in D_{0}^{'}$, it follows that
  $\Op{U}^{'}\tilde\Op{T}_{ij}=\tilde\Op{T}_{ij}\Op{U}^{'}$ and
  $\tilde\Op{T}_{ij}$ is affiliated with $\alg{R}$.

  Further
  \begin{equation*}
    \begin{split}
      \Op{E}_{mn}\Op{S}_{ij}(\sum_{k}\Op{M}_{k}^{'}v_{k})
      &=\sum_{k}\Op{E}_{mn}\Op{E}_{ki}\Op{T}\Op{E}_{jk}\Op{M}_{k}^{'}v_{k}\\
      &=\Op{E}_{mi}\Op{T}\Op{E}_{jm}\Op{E}_{mn}\Op{M}_{n}^{'}v_{n}\\
      &=\sum_{k}\Op{E}_{ki}\Op{T}\Op{E}_{jk}\Op{E}_{mn}\Op{M}_{n}^{'}v_{n}\\
      &=\Op{S}_{ij}\Op{E}_{mn}(\sum_{k}\Op{M}_{k}^{'}v_{k}),
    \end{split}
  \end{equation*}
  hence $\tilde\Op{T}_{ij}$ is affiliated with
  $\alg{T}\otimes\C\otimes\C=
  \{\Op{E}_{mn}\vert m,n\in\N\}^{'}\cap\alg{R}$. Now set
  $\Op{T}_{ij}:=\Op{V}^{*}\tilde\Op{T}_{ij}\Op{V}$, where 
  \begin{equation*}
    \begin{split}
      \Op{V}:\H&\to\raum{K}\\
             v&\mapsto(\delta_{1i}\delta_{i}^{j}v)_{i}^{j}
     \end{split}
  \end{equation*}
  is the canonical partial isometry from $\H$ to
  $\raum{K}=\H\otimes\H_{\infty}\otimes\H_{\infty}$. 

  Now $u_{\tr}\in\Dom{T_{ij}}$ since $\Op{V}u_{\tr}=v_{1}\in\Dom{T_{ij}}$, and
  with $u_{i}^{j}:=\Op{T}_{ij}u_{\tr}$
  \begin{equation*}
    \begin{split}
      \sum_{i,j}\norm{u_{i}^{j}}^{2}
      &=\sum_{i,j}\norm{\Op{T}_{ij}u_{\tr}}^{2}\\
      &=\sum_{i,j}\norm{\Op{V}^{*}\tilde\Op{T}_{ij}\Op{V}u_{\tr}}^{2}\\
      &=\sum_{i,j}\norm{\Op{E}_{1i}\Op{T}\Op{E}_{j1}v_{1}}^{2}\\
      &=\sum_{i,j}\norm{\Op{E}_{i}\Op{T}v_{j}}^{2}\\
      &=\sum_{j}\norm{\Op{T}v_{j}}^{2}<\infty
    \end{split}
  \end{equation*}
  s.t. $u:=\sum_{k}\Op{T}v_{k}=(u_{i}^{j})_{i}^{j}=
           (\Op{T}_{ij}u_{\tr})_{i}^{j}\in\raum{K}$. This means that we can
  construct the operator $\Op{T}_{u}=(\Op{T}_{ij})_{ij}$ according to
  \lemref{lem2:2:2a}. Now
  $\Op{T}_{u}$ and $\Op{T}$ coincide on the core $D_{0}^{'}$, and hence they
  are equal.
\end{proof}

\begin{cor}\label{cor2:2:2}
  The operator $\Op{T}_{u}$ defined in \defnref{defn2:2:1} is unique.
\end{cor}

\begin{cor}\label{cor2:2:3}
  Let $\Op{T}_{u}$ be the operator defined in \defnref{defn2:2:1}. Then
  $\alg{R}v_{k}\in\Dom{T_{u}}$ for every $k\in\N$ and
  \begin{equation*}
    \Op{T}_{u}\Op{M}v_{k}=
    (\sum_{l}\Op{T}_{il}\Op{M}_{lk}\delta_{k}^{j}u_{\tr})_{i}^{j}.
  \end{equation*}
\end{cor}

\begin{proof}
  \propref{prop2:2:2} shows that $\alg{R}v_{k}\in\Dom{T_{u}}$ for every
  $k\in\N$ and $\Op{M}v_{k}=\sum_{l}\Op{M}^{'}_{(k,l)}v_{l}$
  (cf. \eqref{eq2:2:4}). Now
  \begin{equation*}
    \begin{split}
      \Op{T}_{u}\Op{M}v_{k}
      &=\Op{T}_{u}\sum_{l}\Op{M}^{'}_{(k,l)}v_{l}\\
      &=\sum_{l}\Op{T}_{u}\Op{M}^{'}_{(k,l)}v_{l}\\
      &=\sum_{l}(\Op{T}_{il}{\Op{M}^{'}}^{k,l}\delta_{jk}u_{\tr})_{i}^{j}\\
      &=\sum_{l}(\Op{T}_{il}{\Op{M}}_{lk}\delta_{k}^{j}u_{\tr})_{i}^{j}.
    \end{split}
  \end{equation*}
\end{proof}

Now we can formulate the following lemma: 
\begin{lem}\label{lem2:2:2}
  Let $\Op{T}_{u}$ be the operator defined in \defnref{defn2:2:1}. Then:
  \begin{enumerate}
    \item $\tr(\Op{T}_{u}^{*}\Op{T}_{u})=\tr(\Op{T}_{u}\Op{T}_{u}^{*})<\infty$.
    \item $u$ is cyclic, iff $\Op{T}_{u}$ is injective.
    \item $u$ is separating, iff $\Op{T}_{u}$ has dense range.
    \item $u$ is cyclic and separating iff $\Op{T}_{u}$ is injective and has
  dense range, i.e. iff $\Op{T}_{u}$ is invertible. 
  \end{enumerate}
\end{lem}

For the proof we need:
\begin{prop}\label{prop2:2:0}
  Let $\alg{T}$ be a (finite) von Neumann algebra with cyclic trace vector
  $u_{\tr}$. Let further $\Op{S},\Op{T}\eta\alg{T}$ with
  $u_{\tr}\in\Dom{S}\cap\Dom{T}$ and $\Op{M},\Op{N}\in\alg{T}$. Then 
  \begin{equation}\label{eq2:2:3}
    \SProd{\Op{M}\Op{T}u_{\tr}}{\Op{N}\Op{S}u_{\tr}}=
    \SProd{\Op{S}^{*}\Op{N}^{*}u_{\tr}}{\Op{T}^{*}\Op{M}^{*}u_{\tr}}.
  \end{equation}
\end{prop}
\begin{proof}
  Let $(\Op{E}_{n})$ and $(\Op{F}_{n})$ be bounding sequences for $\Op{T}$ and
  $\Op{S}$, resp. (cf. \cite[Lem. 5.6.14]{KRI}). Then:
  \begin{equation*}
    \begin{split}
      \SProd{\Op{M}\Op{T}u_{\tr}}{\Op{N}\Op{S}u_{\tr}}
      &=\Lim{n}\SProd{\Op{M}\Op{T}\Op{E}_{n}u_{\tr}}
               {\Op{N}\Op{S}\Op{F}_{n}u_{\tr}}\\
      &=\Lim{n}\SProd{(\Op{S}\Op{F}_{n})^{*}\Op{N}^{*}u_{\tr}}
                     {(\Op{T}\Op{E}_{n})^{*}\Op{M}^{*}u_{\tr}}\\
      &=\Lim{n}\SProd{\Op{F}_{n}\Op{S}^{*}\Op{N}^{*}u_{\tr}}
                     {\Op{E}_{n}\Op{T}^{*}\Op{M}^{*}u_{\tr}}\\
      &=\SProd{\Op{S}^{*}\Op{N}^{*}u_{\tr}}{\Op{T}^{*}\Op{M}^{*}u_{\tr}},  
    \end{split}
  \end{equation*}
  since $\Op{N}^{*}u_{\tr}\in\Dom{S^{*}}$ and
  $\Op{M}^{*}u_{\tr}\in\Dom{T^{*}}$ (cf. \cite[Prop2.1]{BolI}).
\end{proof}

\begin{proof}[Proof of \lemref{lem2:2:2}]
  \begin{enumerate}
      \item Since $v_{k}\in\Dom{T_{u}}=\Dom{H}$ for all $k\in\N$ , where
        $\Op{T}_{u}=\Op{V}\Op{H}$ is the polar decomposition of $\Op{T}_{u}$,
        we can write the
        trace, defined in \eqref{eq2:2:1}, as follows ($\Op{E}_{\lambda}$ is
        the spectral measure of $\Op{H}$): 
        \begin{equation*}
          \tr(\Op{T}_{u}^{*}\Op{T}_{u})=\tr(\Op{H}^{2})=
          \sum_{k}\int\lambda^{2}d\norm{\Op{E}_{\lambda}v_{k}}^{2}=         
          \sum_{k}\norm{\Op{H}v_{k}}^{2}
          =\sum_{k}\norm{\Op{T}_{u}v_{k}}^{2}.
        \end{equation*}
        Since the $[\alg{R}v_{k}]$ are mutually orthogonal, we have
        \begin{equation*}
          \begin{split}
            \tr(\Op{T}_{u}^{*}\Op{T}_{u})
            &=\sum_{k}\SProd{\Op{T}_{u}v_{k}}{\Op{T}_{u}v_{k}}\\
            &=\sum_{k,j}\SProd{\Op{T}_{u}v_{k}}{\Op{T}_{u}v_{j}}\\
            &=\norm{\sum_{k}\Op{T}_{u}v_{k}}^{2}=\norm{u}^{2}<\infty.
          \end{split}
        \end{equation*}
        Further
        \begin{equation*}
          \begin{split}
            \tr(\Op{T}_{u}\Op{T}_{u}^{*})
            &=\sum_{j}\SProd{\Op{T}_{u}^{*}v_{j}}{\Op{T}_{u}^{*}v_{j}}\\
            &=\sum_{j}\sum_{k,i}
              \SProd{\Op{T}_{jk}^{*}\delta^{ji}u_{\tr}}
                    {\Op{T}_{jk}^{*}\delta^{ji}u_{\tr}}\\
            &=\sum_{k,i}
              \SProd{\Op{T}_{ik}^{*}u_{\tr}}{\Op{T}_{ik}^{*}u_{\tr}}\\
            &=\sum_{k}\sum_{i}
              \SProd{\Op{T}_{ik}u_{\tr}}{\Op{T}_{ik}u_{\tr}}\\
            &=\sum_{k}\SProd{\Op{T}_{u}v_{k}}{\Op{T}_{u}v_{k}}\\
            &=\tr(\Op{T}_{u}^{*}\Op{T}_{u}).
          \end{split}
        \end{equation*}
      \item Let $u$ be cyclic. Then there are 
        $\Op{M}^{(n)}=(\Op{M}_{ik}^{(n)})\in\alg{R}$ with
        \begin{equation*}
          \Lim{n}\Op{M}^{(n)}u=v
        \end{equation*}
        for every
        $v=(\Op{S}_{ij}u_{\tr})_{i}^{j}\in\raum{K}$, where
        $\Op{S}=(\Op{S}_{ij})\in\alg{R}_{0}$. This
        means, using \propref{prop2:2:0} and \corref{cor2:2:3}, 
        \begin{equation*}
          \begin{split}
            0\xleftarrow{\infty\gets n}&
            \sum_{i,j}\norm{\sum_{k}\Op{M}_{ik}^{(n)}\Op{T}_{kj}u_{\tr}-
                  \Op{S}_{ij}u_{\tr}}^{2}\\
            =&\sum_{i,j}(\norm{\sum_{k}\Op{M}_{ik}^{(n)}\Op{T}_{kj}u_{\tr}}^{2}-
              2\sum_{k}\RE\SProd{\Op{M}_{ik}^{(n)}\Op{T}_{kj}u_{\tr}}
                                {\Op{S}_{ij}u_{\tr}}+
              \norm{\Op{S}_{ij}u_{\tr}}^{2})\\
            =&\sum_{i,j}(\norm{\sum_{k}\Op{T}_{kj}^{*}(\Op{M}_{ik}^{(n)})^{*}u_{\tr}}^{2}-
              2\sum_{k}\RE\SProd{\Op{T}_{kj}^{*}(\Op{M}_{ik}^{(n)})^{*}u_{\tr}}
                                {\Op{S}_{ij}^{*}u_{\tr}}+
              \norm{\Op{S}_{ij}^{*}u_{\tr}}^{2})\\
            =&\sum_{i,j}\norm{\sum_{k}\Op{T}_{kj}^{*}(\Op{M}_{ik}^{(n)})^{*}u_{\tr}-
                  \Op{S}_{ij}^{*}u_{\tr}}^{2},
          \end{split}
        \end{equation*}
        i.e., since $(\Op{T}_{ki}^{*})_{i,k}\subset\Op{T}_{u}^{*}$ and
        $\alg{R}_{0}u_{\tr}$ is dense in $\raum{K}$, $\Op{T}_{u}^{*}$ has
        dense range, i.e. $\Op{T}_{u}$ is injective.

        Let now $\Op{T}_{u}$ be injective and
        $\Op{M}^{'}=({\Op{M}^{'}}^{ij})\in\alg{R}^{'}$ with
        $\Op{M}^{'}u=0$. Now 
        \begin{equation*}
          \Op{M}^{'}u=(\sum_{j}{\Op{M}^{'}}^{ij}\Op{T}_{kj}u_{\tr})_{ik}=0,
        \end{equation*}
        and
        \begin{equation*}
          \begin{split}
            0=\norm{\Op{M}^{'}u}^{2}
            &=\sum_{i,k}\norm{\sum_{j}{\Op{M}^{'}}^{ij}\Op{T}_{kj}u_{\tr}}\\
            &=\sum_{i,k}\norm{\sum_{j}\Op{T}_{kj}{\Op{M}^{'}}^{ij}u_{\tr}}
            =\norm{\Op{T}_{u}v},
          \end{split}
        \end{equation*}
        where
        $v:=({\Op{M}^{'}}^{ij}u_{\tr})_{i}^{j}
           =\sum_{k}\Op{M}^{'}v_{k}\in\Dom{T_{u}}$ ($\Op{T}_{u}$ is closed),
        hence $\Op{T}_{u}v=0$, and, since $\Op{T}_{u}$ is injective,
        $v_{i}^{j}={\Op{M}^{'}}^{ij}u_{\tr}=0$ for all $i,j\in\N$. Because
        $u_{\tr}$ is cyclic for $\alg{T}$ hence separating for $\alg{T}^{'}$,
        ${\Op{M}^{'}}^{ij}=0$ for all $i,j\in\N$, s.t. $\Op{M}^{'}=0$.
      \item Let $u$ be separating. This means that $u$ is cyclic for
        $\alg{R}^{'}$. Then there are 
        $\Op{M}_{(n)}=(\Op{M}_{(n)}^{ik})\in\alg{R}^{'}$ and
        \begin{equation*}
          \Lim{n}\Op{M}_{(n)}u=v
        \end{equation*}
        for every $v=(\Op{S}_{ij}u_{\tr})_{i}^{j}\in\raum{K}$,
        where $(\Op{S}_{ij})\in\alg{R}_{0}^{'}$ ($i,j\in\N$). This means
        \begin{equation*}
          \begin{split}
            0\xleftarrow{\infty\gets n}&
            \sum_{i,j}\norm{\sum_{k}\Op{M}_{(n)}^{jk}\Op{T}_{ik}u_{\tr}-
                  \Op{S}_{ij}u_{\tr}}^{2}\\
            =&\sum_{i,j}\norm{\sum_{k}\Op{T}_{ik}\Op{M}_{(n)}^{jk}u_{\tr}-
                  \Op{S}_{ij}u_{\tr}}^{2}.
          \end{split}
        \end{equation*}
        Since $\alg{R}_{0}^{'}u_{\tr}$ is dense in $\raum{K}$ we have
        proven that $\Op{T}_{u}$ has dense range.

        For the converse read the argument backwards.

      \item This follows from 2. and 3.
    \end{enumerate}
\end{proof}

\begin{rem}\label{rem2:2:1}
  Also here, as in the finite case, the finite trace
  condition of \lemref{lem2:2:2} is not only necessary but also sufficient
  for an operator being the operator associated with a vector in the sense of
  \defnref{defn2:2:1}. Suppose that $\tr(\Op{T}^{*}\Op{T})<\infty$ with
  $\Op{T}\eta\alg{R}$. Then
  \begin{equation*}
    \begin{split}
      \infty>\tr(\Op{T}^{*}\Op{T})
      &=\tr(\Op{H}^{2})\\
      &=\sum_{k}\int\lambda^{2}d\norm{\Op{E}_{\lambda}v_{k}}
    \end{split}
  \end{equation*}
  hence
  \begin{equation*}
    \int\lambda^{2}d\norm{\Op{E}_{\lambda}v_{k}}<\infty\quad\forall k\in\N,
  \end{equation*}
  i.e. $v_{k}\in\Dom{H}=\Dom{T}$, and
  \begin{equation*}
    \sum_{k}\norm{\Op{T}v_{k}}^{2}=\sum_{k}\norm{\Op{H}v_{k}}^{2}
    =\sum_{k}\int\lambda^{2}d\norm{\Op{E}_{\lambda}v_{k}}<\infty.
  \end{equation*}
  This shows that the assumptions of \lemref{lem2:2:3} are fulfilled.
\end{rem}

\begin{cor}\label{cor2:2:1}
  $\alg{R}$ possesses a cyclic and separating vector $u_{0}\in\raum{K}$.
\end{cor}

\begin{proof}
  Set $\Op{T}:=(\delta_{ij}j^{-2}\Op{Id})_{i,j}$ or
  $u_{0}:=\sum_{j}j^{-2}v_{j}$. Then $\Op{T}$ fulfills the conditions of
  \lemref{lem2:2:3} and is invertible, s.t. from \lemref{lem2:2:2} follows
  that $u_{0}$ is cyclic and separating.
\end{proof}

In the last step of this subsection we show that the model we have just
treated is really representative for the general situation, in the sense that
all infinite type $I$ or type $II$ factors can be considered as a matrix
algebra of finite type $I$ or type $II$ factors, resp. This is shown by the
next 
\begin{lem}\label{lem2:2:1}
  Every infinite but semifinite von Neumann factor $(\M_{0},\H_{0})$ with
  cyclic and separating 
  vector $u_{0}\in\H_{0}$ is unitarily equivalent to
  $\alg{T}\otimes L(\H_{\infty})\otimes\C=:\alg{R},
    \H\otimes\H_{\infty}\otimes\H_{\infty}=:\raum{K})$,
  where $\alg{T}$ is a finite von Neumann factor acting on the Hilbert space
  $\H$ with cyclic and separating vector and
  $\H_{\infty}$ is a separable infinite dimensional Hilbert space. 
\end{lem}

\begin{proof}
  Since $\M_{0}$ is infinite but semifinite there is a countable orthogonal
  family of finite equivalent projections $(\Op{E}_{n})_{n\in\N}$ in $\M_{0}$,
  s.t. $\sum\Op{E}_{n}=\Op{Id}$. Now there is a selfadjoint system of matrix
  units $(\Op{E}_{ab})_{a,b\in\N}$ with $\Op{E}_{aa}=\Op{E}_{a}$
  (cf. \cite[6.6.4]{KRII}). This shows that $\M_{0}$ is isomorphic to
  $\tilde\alg{T}\otimes L(\H_{\infty})$ where
  $\tilde\alg{T}:=\{\Op{E}_{ab}\}^{'}\cap\alg{M}_{0}$ and $\tilde\alg{T}$ is
  isomorphic to every $\Op{E}_{n}\M_{0}\Op{E}_{n}$ ($n\in\N$). Since the
  projections $\Op{E}_{n}$ are finite also $\tilde\alg{T}$ is a finite factor.

  Since $\M_{0}$ possesses the separating vector $u_{0}$ we can represent the
  algebras $\Op{E}_{n}\M_{0}\Op{E}_{n}$ by the GNS representation for the
  faithful state $\omega_{n}$ induced by the
  separating vector $\Op{E}_{n}u_{0}$ on a Hilbert space $\H_{n}$,
  s.t. the vector $u_{n}\in\H_{n}$ implementing the state $\omega_{n}$ is a
  cyclic and separating vector for $\Op{E}_{n}\M_{0}\Op{E}_{n}$. Since all the
  $\Op{E}_{n}\M_{0}\Op{E}_{n}$ are isomorphic and they possess in this
  representation a cyclic and separating vector, they are all unitarily
  equivalent. This means that we can choose as  $\alg{T}$ one of the
  $\Op{E}_{n}\M_{0}\Op{E}_{n}$ acting on the representation space $\H_{n}$.

  Since the factor $(\alg{T}\otimes L(\H_{\infty})\otimes\C=:\alg{R},
    \H\otimes\H_{\infty}\otimes\H_{\infty}=:\raum{K})$ possesses a cyclic and
    separating vector if $(\alg{T},\H)$ does (see \corref{cor2:2:1}) and it
    is isomorphic to $\M_{0}$ it is unitarily equivalent to $\M_{0}$.
\end{proof}

The results of this section (and the analogues in \cite{BolI}) can be subsumed
in the next two theorems:
\begin{thm}\label{thm2:1}
  Let ($\M_{0},\H_{0}$) be a semifinite von Neumann factor. Let further
  $u\in\H_{0}$. Then there is
  exactly one operator $\Op{T}_{u}\eta\M_{0}$ associated with the vector $u$
  in the sense of \cite[Def 2.1.]{BolI} in the finite case and in the sense of
  \defnref{defn2:2:1} in the infinite case, resp., having the
  following properties: 
  \begin{enumerate}
    \item $\tr(\Op{T}_{u}\Op{T}_{u}^{*})=\tr(\Op{T}_{u}^{*}\Op{T}_{u})<\infty$.
    \item $u$ is cyclic, iff $\Op{T}_{u}$ is injective.
    \item $u$ is separating, iff $\Op{T}_{u}$ has dense range.
    \item $u$ is cyclic and separating iff $\Op{T}_{u}$ is injective and has
  dense range, i.e. iff $\Op{T}_{u}$ is invertible. 
  \end{enumerate}
\end{thm}
\begin{proof}
  The finite case is just \thmref{thm1:1}. In the infinite case the existence
  and the asserted properties follow from \lemref{lem2:2:1} and
  \lemref{lem2:2:2} infinite case, the uniqueness from \corref{cor2:2:2}.
\end{proof}

\begin{thm}\label{thm2:2}
  Let $\Op{T}\eta\M_{0}$. Then there is a vector $u\in\H_{0}$
  s.t $\Op{T}=\Op{T}_{u}$ iff
  $\tr(\Op{T}\Op{T}^{*})=\tr(\Op{T}^{*}\Op{T})<\infty$.  
\end{thm}

\begin{proof}[Proof of \thmref{thm2:2}]
  Again the finite case is just \thmref{thm1:2}. In the infinite case the
  necessarity of the trace condition follows from \thmref{thm2:1} and the
  sufficiency from \remref{rem2:2:1}, resp.
\end{proof}

\section{Generation of Modular Objects}
\label{sec:3}

In this section we show how the modular objects of a cyclic and separating
vector $u_{0}\in\H$ for a semifinite von Neumann factor $(\M_{0},\H_{0})$ are
related to the operator $\Op{T}_{u_{0}}$ constructed in the last section. As
in \secref{sec:2} we consider as a model for the infinite but semifinite 
factor the factor $\alg{T}\otimes L(\H_{\infty})\otimes\C=:\alg{R},
    \H\otimes\H_{\infty}\otimes\H_{\infty})=:\raum{K})$, where $\alg{T}$ is a
finite factor with cyclic trace vector $u_{\tr}\in\H$. If $u_{0}\in\raum{K}$
is a cyclic and separating vector for $\alg{R}$, according to
\lemref{lem2:2:2}, there is an invertible operator $\Op{T}_{u_{0}}\eta\alg{R}$,
s.t. $u_{0}=\sum_{k}\Op{T}_{u_{0}}v_{k}$, where
$v_{k}=(\delta_{i}^{j}\delta_{ik}u_{\tr})_{i}^{j}$. Using this operator we can
formulate the following analogue to \thmref{thm1:3}:
\begin{thm}\label{thm3:1}
  Use the notations from above. Let further
  \begin{equation*}
    \Op{T}_{u_{0}}=\Op{HV}=
    (\Op{H_{ij}})_{ij}
    (\Op{V}_{ij})_{ij}
  \end{equation*}
  be the polar decomposition of $\Op{T}_{u_{0}}$. With the conjugation
  $\tilde\Op{J}$ defined as
  \begin{equation*}
    \tilde\Op{J}(\Op{M}_{ij}u_{\tr})_{i}^{j}
    :=(\Op{M}_{ji}^{*}u_{\tr})_{i}^{j}
    :=(\Op{J}\Op{M}_{ji}u_{\tr})_{i}^{j}
    \quad\forall\Op{M}=(\Op{M}_{ij})_{ij}\in\alg{R},
  \end{equation*}
  where $\Op{J}$ is the conjugation corresponding to the trace vector
  $u_{\tr}$, we can calculate the
  modular objects $(\Delta_{0},\Op{J}_{0})$ of $(\M_{0},u_{0})$ as follows:
  \begin{equation*}
    \Op{J}_0=
    \tilde\Op{J}\Op{V}^*\tilde\Op{J}\Op{V}\tilde\Op{J}=
    \Op{V}\tilde\Op{J}\Op{V}^*,
  \end{equation*}
  and
  \begin{equation*}
    \Delta_0=\Op{J}_0\Op{H}^{-1}_0\Op{J}_{0}\Op{H}_0,
  \end{equation*}
  where $\Op{H}_0=\Op{H}^2=\Op{T}_{u_{0}}\Op{T}_{u_{0}}^{*}$.
\end{thm}

\begin{proof}
  \begin{enumerate}
    \item
      First we observe that $\tilde\Op{J}\Op{R}\tilde\Op{J}\in\alg{R}^{'}$ for
      every $\Op{R}\in\alg{R}$. For let
      $\Op{R}=(\Op{R}_{ij})\in\alg{R}$ and
      $v=(u_{i}^{j})_{i}^{j}=(\Op{M}_{ij}u_{\tr})_{i}^{j}\in\raum{K}$,
      $(\Op{M}_{ij})\in\alg{R}_{0}$, then
      \begin{equation*}
        \begin{split}
          \tilde\Op{J}\Op{R}\tilde\Op{J}v
          &=\tilde\Op{J}\Op{R}(\Op{J}\Op{M}_{ji}u_{\tr})_{i}^{j}\\
          &=\tilde\Op{J}(\sum_{i}\Op{R}_{ki}\Op{J}\Op{M}_{ji}u_{\tr})_{k}^{j}\\
          &=(\Op{J}\sum_{i}\Op{R}_{ji}\Op{J}\Op{M}_{ki}u_{\tr})_{k}^{j}\\
          &=\underbrace{(\Op{J}\Op{R}_{ji}\Op{J})^{ji}}
          _{:=\Op{R}^{'}\in\alg{R}^{'}}
          (\Op{M}_{ki}u_{\tr})_{k}^{i}\\
          &=\Op{R}^{'}v.
        \end{split}
      \end{equation*}
      Further
      \begin{equation*}
        \tilde\Op{J}\tilde\Op{J}v
        =\tilde\Op{J}(\Op{J}\Op{M}_{ji}u_{\tr})_{ij}
        =(\Op{M}_{ij}u_{\tr})_{ij}=v,
      \end{equation*}
      s.t. $\tilde\Op{J}$ is an (algebraic) conjugation for $\alg{R}$.

    \item Let now $\Op{T}_{u_{0}}$ be bounded ($\Rightarrow$ all the
      $\Op{T_{ij}}$ and $\Op{H}_{ij}$, resp. are bounded). Then we show that
      the Tomita operator $\Op{S}$ defined by
      \begin{equation*}
        \Op{S}\Op{A}u_{0}=\Op{A}^{*}u_{0}\quad\forall\Op{A}\in\alg{R}
      \end{equation*}
      can be written as
      \begin{equation}\label{eq3:2:1}
        \Op{S}=\Op{H}^{-1}\Op{V}\tilde\Op{J}\Op{V}^{*}\Op{H}.
      \end{equation}
      For this let $\Op{A}=(\Op{A}_{ij})_{ij}\in\alg{R}$
      and $u_{0}=(\sum_{k}\Op{H}_{jk}\Op{V}_{kl}u_{\tr})_{j}^{l}$. Then
      \begin{equation*}
        \Op{A}u_{0}=
        (\sum_{j,k}\Op{A}_{ij}\Op{H}_{jk}\Op{V}_{kl}u_{\tr})_{i}^{l}
      \end{equation*}
      and
      \begin{equation*}
        \Op{A}^{*}u_{0}=
        (\sum_{j,k}\Op{A}_{ji}^{*}\Op{H}_{jk}\Op{V}_{kl}u_{\tr})_{i}^{l}.
      \end{equation*}
      Now
      \begin{equation*}
        \begin{split}
          (\Op{H}^{-1}\Op{V}\tilde\Op{J}\Op{V}^{*}\Op{H})\Op{A}u_{0}
          &=\Op{H}^{-1}\Op{V}\tilde\Op{J}
            (\sum_{i,j,k,m}
            \Op{V}_{mn}^{*}\Op{H}_{mi}\Op{A}_{ij}
            \Op{H}_{jk}\Op{V}_{kl}u_{\tr})_{n}^{l}\\
          &=\Op{H}^{-1}\Op{V}
            (\sum_{i,j,k,m}\Op{J}
            \Op{V}_{ml}^{*}\Op{H}_{mi}\Op{A}_{ij}
            \Op{H}_{jk}\Op{V}_{kn}u_{\tr})_{n}^{l}\\
          &=\Op{H}^{-1}\Op{V}
            (\sum_{i,j,k,m}
            \Op{V}_{kn}^{*}\Op{H}_{jk}^{*}\Op{A}_{ij}^{*}
            \Op{H}_{mi}^{*}\Op{V}_{ml}u_{\tr})_{n}^{l}\\
          &=\Op{H}^{-1}\Op{V}
            (\sum_{i,j,k,m}
            \Op{V}_{kn}^{*}\Op{H}_{kj}\Op{A}_{ij}^{*}
            \Op{H}_{im}\Op{V}_{ml}u_{\tr})_{n}^{l}\\
          &=\Op{H}^{-1}
            (\sum_{i,j,m}
            \Op{H}_{nj}\Op{A}_{ij}^{*}
            \Op{H}_{im}\Op{V}_{ml}u_{\tr})_{n}^{l}\\
          &=(\sum_{i,m}\Op{A}_{in}^{*}
            \Op{H}_{im}\Op{V}_{ml}u_{\tr})_{n}^{l}=\Op{A}^{*}u_{0},
        \end{split}
      \end{equation*}
      which proves \eqref{eq3:2:1}. Now
      $\Op{S}^{*}=\Op{H}\Op{V}\tilde\Op{J}\Op{V}^{*}\Op{H}^{-1}$ and
      \begin{equation*}
        \begin{split}
          \Delta_{0}
          &=\Op{S}^{*}\Op{S}\\
          &=\Op{H}\Op{V}\tilde\Op{J}\Op{V}^{*}\Op{H}^{-1}
            \Op{H}^{-1}\Op{V}\tilde\Op{J}\Op{V}^{*}\Op{H}\\
          &=\Op{V}\tilde\Op{J}\Op{V}^{*}\Op{H}^{-2}
            \Op{V}\tilde\Op{J}\Op{V}^{*}\Op{H}^{2}\\
          &=\Op{J}_{0}\Op{H}_{0}^{-1}\Op{J}_{0}\Op{H}_{0}.
        \end{split}
      \end{equation*}
      Further 
      \begin{equation*}
        \Op{J}_{0}\Delta_{0}^{1/2}=\Op{H}^{-1}\Op{J}_{0}\Op{H}=\Op{S},
      \end{equation*}
      and all the assertions are proven in the bounded case.
    \item
      In the last step we approximate the (unbounded) operator
      $\Op{T}_{u_{0}}$ by bounded operators $\Op{T}_{n}$ in exactly the same
      way as in the proof of Theorem 3.1. in \cite{BolI} and show the
      assertions like there also in the unbounded case.
  \end{enumerate}
\end{proof}

\section{The Second Simple Class of Solutions of the Inverse Problem}
\label{sec:4}

In this section we want to use the results of the last two sections to examine
the second simple classes of solutions of the inverse problem introduced by
Wollenberg in \cite{WolII} for type $I$ factors, and considered in \cite{BolI}
also for type $II_{1}$ factors. For the construction of this class it is
crucial that the inverse $\Delta_{0}^{-1}$ of the modular operator is again a
modular operator. To this scope there was shown the following
\begin{lem}\label{lem4:1}
  Let $\Delta_{0}=\Op{J}_{0}\Op{H}_{0}^{-1}\Op{J}_{0}\Op{H}_{0}$ be the
  decomposition of the modular operator $\Delta_{0}$, where
  $\Op{J}_0=\Op{JV}^*\Op{JVJ}=\Op{VJV}^*$ and
  $\Op{T}_{u_{0}}=\Op{H}_{0}^{1/2}\Op{V}$ is the operator corresponding to 
  $u_{0}$ (cf. \thmref{thm3:1}). Then
  $\Delta_{0}^{-1}=\Op{J}_{0}\Op{H}_{0}\Op{J}_{0}\Op{H}_{0}^{-1}$ and the
  following is equivalent:
  \begin{enumerate}
  \item ($\Delta_{0}^{-1},\Op{J}_{0})$ are the modular objects w.r.t. a
    cyclic and separating vector $u_{1}\in\H_{0}$.
  \item  
    \begin{equation}\label{eq4:3}
      \tr(\Op{H}_{0}^{-1})<\infty.
    \end{equation}
  \end{enumerate}
\end{lem}
This lemma can be proven with the same techniques as in \cite{BolI} also for
the infinite case taking into account \thmref{thm2:1}, \thmref{thm2:2}, and
\thmref{thm3:1}.

Now we must examine, whether or not the second condition in \lemref{lem4:1} is
fulfilled:
\begin{lem}\label{lem4:2}
  For type $I_{\infty}$ and type $II_{\infty}$ factors the condition
  \eqref{eq4:3} is never true.
\end{lem}

\begin{proof}
  Let $\M_{0}$ now be a type $I_{\infty}$ or $II_{\infty}$ factor and
  $\Op{T}_{u_{0}}=\Op{H}_{0}^{1/2}\Op{V}$ the operator corresponding to
  the cyclic and separating vector $u_{0}$. Let further
  $\Op{E}_{\lambda}\in\M_{0}$ the spectral resolution of $\Op{H}_{0}$. Then
  we can define a positive measure $\mu_{\tr}$ on the $\sigma$-algebra of
  Borel sets in $\R$, s.t. 
  \begin{equation*}
    \tr(\Op{H}_{0})=\int\lambda
    d\mu_{\tr}(\lambda),
  \end{equation*}
  where
  \begin{equation*}
    \mu_{\tr}(B):=\tr{\Op{E}(B)}
  \end{equation*}
  for all Borel sets $B$. Now $c:=\tr(\Op{H}_{0})<\infty$. Assume
  w.l.o.g. $c=1$. Then
  \begin{equation*}
    1=\int\lambda d\mu_{\tr}(\lambda)
    \geq\int_{\abge{0}{1}}\lambda d\mu_{\tr}(\lambda)+
    \int_{\offen{1}{\infty}}d\mu_{\tr}(\lambda),
  \end{equation*}
  i.e.
  \begin{equation*}
    \int_{\offen{1}{\infty}}d\mu_{\tr}(\lambda)<\infty.
  \end{equation*}
  Since $\M_{0}$ is infinite $\infty=\tr(\Op{Id})=\mu_{\tr}(\R)$, i.e.
  \begin{equation*}
    \infty=\int_{\lambda}d\mu_{\tr}(\lambda)=
    \int_{\abge{0}{1}}d\mu_{\tr}(\lambda)+
    \underbrace{\int_{\offen{1}{\infty}}d\mu_{\tr}(\lambda)}_{<\infty},
  \end{equation*}
  hence
  \begin{equation*}
    \int_{\abge{0}{1}}d\mu_{\tr}(\lambda)=\infty.
  \end{equation*}
  Suppose now that also $\tr(\Op{H}_{0}^{-1})<\infty$, then
  \begin{equation*}
    \infty>\int_{\lambda}\lambda^{-1}d\mu_{\tr}(\lambda)
    \geq
    \underbrace{\int_{\abge{0}{1}}d\mu_{\tr}(\lambda)}_{=\infty}+
    \int_{\offen{1}{\infty}}\lambda^{-1}d\mu_{\tr}(\lambda),
  \end{equation*}
  which is a contradiction.
\end{proof}

Hence the last lemma shows that for infinite semifinite factors the second
class of solutions of the inverse problem can never be constructed. This
result was yet obtained by Wollenberg in \cite{WolII} for the type
$I_{\infty}$ case.

\section{The Classification of Solutions in the Pure Point Spectrum Case}
\label{sec:5}

In this section we want to show the modifications of classification results
obtained in \cite{BolI}. The definition of the equivalence relation does not
use any special properties of the finite factors, and can just be repeated
here: 
\begin{defn}\label{defn5:1}
  Two semifinite von Neumann factors
  $\M,\alg{N}\in NF_{\M_{0}}(\Delta_{0},\Op{J}_{0},u_{0})$ are called
  equivalent, $\M\sim\alg{N}$, if $\M\in
  NF^{1}_{\alg{N}}(\Delta_{0},\Op{J}_{0},u_{0})$, i.e. if there exists
  a unitary operator $\Op{U}$ on $\H_{0}$,
  s.t. $\M=\Op{U}\alg{N}\Op{U}^{*}$, $\Op{U}$ commutes with $\Delta_{0}$ and
  $\Op{J}_{0}$ and $\Op{U}^{*}u_{0}=\pm u_{0}$ (For the definition of the
  class $NF^{1}_{\alg{N}}(\Delta_{0},\Op{J}_{0},u_{0})$ see \cite{BolI}).
\end{defn}

Also the next lemmas can be formulated and proved in exactly the same way as
in the finite case. Assume in
the following that $\Op{H}_{0}$ has pure point spectrum, i.e. 
$\Op{H}_{0}=\sum_{k\in K}\mu_{k}\Op{E}_{k}$ where the $\mu_{k}$ ($k\in K$) are
the eigenvalues of
$\Op{H}_{0}$ and $\Op{E}_{k}\in\M_{0}$ are the corresponding (orthogonal)
eigenprojections with
$m_{k}:=\tr{\Op{E}_{k}}=:D_{\M_{0}}(\Op{E}_{k})$ their von Neumann
dimension.
Then we have for $\Delta_{0}$ the following decomposition
\begin{equation}\label{eq5:3}
  \begin{split}
    \Delta_{0}&=\Op{H}_{0}\Op{J}_{0}\Op{H}_{0}^{-1}\Op{J}_{0}\\
              &=\sum_{k,l\in K}\mu_{k}\mu_{l}^{-1}\Op{E}_{k}
                         \Op{J}_{0}\Op{E}_{l}\Op{J}_{0}\\
              &=\sum_{j\in J}\lambda_{j}\Op{F}_{j},
  \end{split}
\end{equation}
where the $\lambda_{j}$ ($j\in J$) are the eigenvalues of $\Delta_{0}$ and
$\Op{F}_{j}$ are the corresponding eigenprojections. Now
\begin{lem}\label{lem5:1} 
  With the notations introduced above we can compute the spectrum of
  $\Delta_{0}$  in the following way: 
  \begin{equation}\label{eq5:4}
    \{\lambda_{j}\vert j\in J\}=\{\mu_{k}\mu_{l}^{-1}\vert k,l\in K\}
    \quad\forall j\in J
  \end{equation}
  and
  \begin{subequations}\label{eq5:5}
    \begin{equation}\label{eq5:5a}
      n_{j}=\sum_{\mu_{k}\mu_{l}^{-1}=\lambda_{j}}m_{k}m_{l}
      \quad\forall j\in J\text{ if $\M_{0}$ is type $I$,}
    \end{equation}
    \begin{equation}\label{eq5:5b}
      n_{j}=\infty\quad\forall j\in J\text{ if $\M_{0}$ is type $II$,}
    \end{equation}  
  \end{subequations}
  where $n_{j}:=D_{L(\H_{0})}(\Op{F}_{j})$ with $D_{L(\H_{0})}(\Op{F}_{j})$
  the dimension function in
  the type $I_{\infty}$ factor $L(\H_{0})$, which corresponds to the
  normalized Hilbert space dimension.
\end{lem}

\begin{lem}\label{lem5:2}
  If there are two solutions of the inverse problem $\M_{1}$, $\M_{2}$ s.t. the
  corresponding selfadjoint operators $\Op{H}_{1}$ and $\Op{H}_{2}$ have the
  same eigenvalues modulo a positive constant $c>0$ and same (von Neumann)
  multiplicities, then $\M_{1}\sim\M_{2}$. 
\end{lem}

\begin{lem}\label{lem5:3}
  If there are two equivalent solutions $\M_{1}$, $\M_{2}$ of the inverse
  problem with the corresponding positive operators $\Op{H}_{1}$ and
  $\Op{H}_{2}$, resp., (having pure point spectrum) then
  $\Op{H}_{1}$ and $\Op{H}_{2}$ have the same eigenvalues (up to a positive
  constant) and von Neumann
  multiplicities, i.e. they are unitarily equivalent in $\M_{0}$.
\end{lem}

The only difference to the finite case is shown by the next  
\begin{lem}\label{lem5:4}
  Let $(\mu_{k},m_{k})_{k\in K}$ be a sequence of pairs of positive reals
    $\mu_{k}>0$ and $m_{k}>0$, s.t.
    \begin{subequations}\label{eq5:6}
      \begin{equation}
        m_{k}\in\N\text{ if $\M_{0}$ is type $I_{\infty}$},
      \end{equation}
      \begin{equation}
        m_{k}\in\R_{>0}\text{ if $\M_{0}$ is type $II_{\infty}$},
      \end{equation}
    and
    \begin{equation}
      \sum_{k\in K}m_{k}=\infty
    \end{equation}
    and 
    \begin{equation}
      \sum_{k\in K}m_{k}\mu_{k}=1
    \end{equation}
  \end{subequations}
  and the relations \eqref{eq5:4} and
  \eqref{eq5:5} are fulfilled. Then there exists a solution
  $\M=\Op{U}\M_{0}\Op{U}^{*}\in NF_{\M_{0}}(\Delta_{0},\Op{J}_{0},u_{0})$,
  s.t.
  $\Op{U}^{*}\Delta_{0}\Op{U}=\Op{H}\Op{J}_{0}\Op{H}^{-1}\Op{J}_{0}$ 
  and $\Op{H}$ has the eigenvalues and multiplicities $(\mu_{k},m_{k})_{k\in
    K}$ (cf. \cite[prop.4.1]{WolII}).
\end{lem}

For the proof we need the following auxiliary results:
\begin{prop}\label{prop5:4}
  If $(m_{k})$ is countable family of positive reals with $\sum m_{k}=\infty$, 
  then there exists in a type $II_{\infty}$ von Neumann factor $\M$ a family of
  pairwise orthogonal projections $(\Op{E}_{k})$, s.t. $D(\Op{E}_{k})=m_{k}$
  for every $k$.
\end{prop}

\begin{proof}
  We construct the $\Op{E}_{k}$ inductively: Since the range of $D_{\M}$ is
  all of $\R_{\le0}$ (cf. \cite[8.4.4]{KRII}) there is a projection in $\M$,
  s.t. $D(\Op{E}_{1})=m_{1}$. 

  Suppose now that for $N\in\N$ the $\Op{E}_{k}$ are pairwise orthogonal with
  $D_{\M_{0}}(\Op{E}_{k})=m_{k}$ ($1\leq k<N$). Setting
  $\Op{F}_{N}:=\Op{Id}-\sum_{k=1}^{N}\Op{E}_{k}$ the restricted algebra
  $\Op{F}_{N}\M\Op{F}_{N}$ is again a type $II$ factor, finite, if
  $\Op{F}_{N}$ is finite, and infinite, if $\Op{F}_{N}$ is infinite (cf.
  \cite[Ex. 6.9.16]{KRII}) with the
  dimension function
  \begin{equation*}
    D_{N}(\Op{F}_{n}\Op{E}\Op{F}_{N}):=
    D_{\M_{0}}(\Op{F}_{n}\Op{E}\Op{F}_{N})/D(\Op{F}_{N})
    \quad\forall\Op{F}_{n}\Op{E}\Op{F}_{N}\in\Op{F}_{N}\M\Op{F}_{N},
  \end{equation*}
  if $\Op{F}_{N}$ is finite, and $D_{N}=D_{\M_{0}}$ else,
  where 
  \begin{equation*}
    D_{\M_{0}}(\Op{F}_{N})=D_{\M_{0}}(\Op{Id}-\sum_{k=1}^{N}\Op{E}_{k})
    =1-\sum_{k=1}^{N}D_{\M_{0}}(\Op{E}_{k})\geq m_{N}.
  \end{equation*}
  With the same argument as above there is again a projection
  $\Op{E}_{N}\in\Op{F}_{N}\M\Op{F}_{N}\subset\M$,
  s.t. $D_{N}(\Op{E}_{N})=D(\Op{F}_{N})^{-1}m_{N}\leq1$, if $\Op{F}_{N}$ is
  finite, and $D_{N}(\Op{E}_{N})=m_{N}$ else. In both cases
  $D_{\M_{0}}(\Op{E}_{N})$ and $\Op{E}_{N}<\Op{F}_{N}\perp\Op{E}_{k}$ ($1\leq
  k<N$).  
\end{proof}

Now the proof of \lemref{lem5:4} is the same as in \cite{BolI}.

\begin{rem}\label{rem5:1}
  \eqref{eq5:6} show that in the infinite case we have always an infinite
    set of eigenvalues with $0$ as cummulation point, i.e. $K=\N$ and $0$ is
    in the spectrum $\sigma(\Op{H})$ of $\Op{H}$.
\end{rem}

Now we can summarize the lemmas of this section in the following
\begin{thm}\label{thm5:1}
  Let $\M_{0}$ be a semifinite von Neumann factor with cyclic and separating
  vector $u_{0}$ and $\Op{T}_{u_{0}}=\Op{H}_{0}^{-1/2}\Op{V}$ the operator
  corresponding to $u_{0}$. If $\Op{H}_{0}$ has pure point spectrum, also
  $\Delta_{0}$ has it. In this case let $(\lambda_{j})$ $(j\in J)$ be the
  eigenvalues of $\Delta_{0}$. Then  
  \begin{enumerate}
  \item Two solutions $\M_{1},M_{2}\in
    NF_{\M_{0}}(\Delta_{0},\Op{J}_{0},u_{0})$ of the inverse problem with
    corresponding invertible operators $\Op{H}_{i}\eta\M_{0}$ $(i=1,2)$ having
    pure point spectrum are equivalent
    iff $\Op{H}_{1}$ and $\Op{H}_{2}$ have the same eigenvalues and (von
    Neumann) multiplicities.
  \item A positive invertible operator $\Op{H}\eta\M_{0}$ with pure point
    spectrum gives rise to a solution of
    the inverse problem iff its eigenvalues and multiplicities satisfy
    \eqref{eq5:4}, \eqref{eq5:5}, and \eqref{eq5:6}.
  \item When the corresponding operators $\Op{H}$ has pure point spectrum the
    equivalence classes of $\sim$ are completely classified by the spectrum of
    the corresponding operators, i.e. by
    sequences of pairs of positive reals $(\mu_{k},m_{k})$ satisfying
    \eqref{eq5:4}, \eqref{eq5:5}, and \eqref{eq5:6}.
  \end{enumerate}
\end{thm}

\begin{bsp}\label{bsp5:1}
  Here we want to give some examples to illustrate \thmref{thm5:1}.
  \begin{enumerate}
  \item In \cite{WolII} you can find some examples for the type $I$ case.
  \item Let 
    \begin{equation*}
      (\ldots,10^{-3},10^{-2},10^{-1},1,10,10^{2},10^{3},\ldots)
    \end{equation*} 
    be the eigenvalues of a modular operator for a type $II_{\infty}$
    factor. Then 
    \begin{equation*}
      ((c_{1}\cdot 1,1),(c_{1}\cdot 10^{-1},1),(c_{1}\cdot
      10^{-2},1),(c_{1} \cdot 10^{-3},1),\ldots)
    \end{equation*}
    and
    \begin{equation*}
      ((c_{2}\cdot 1,1),(c_{2}\cdot 10^{-1},1),(c_{2}\cdot
      10^{-3},1),(c_{2} \cdot 10^{-5},1),\ldots)
    \end{equation*}
    characterize two different classes of solutions of the inverse problem,
    i.e. they both satisfy \eqref{eq5:4}, \eqref{eq5:5}, and \eqref{eq5:6},
    where $c_{i}$ ($i=1,2$) are appropriate chosen constants. This
    shows that in this case there are more than the simple classes of
    solutions of the inverse problem. 
  \item Let $(\mu_{k},m_{k})_{k\in\N}$ characterize a class of solutions of
    the inverse problem in the type $II_{\infty}$ case, where $m_{l}\not=m_{k}$
    for at least one pair $k,l\in\N$, then for every finite permutation
    $\sigma$ of $\N$ interchanging $k$ and $l$ also
    $(c\mu_{k},m_{\sigma(k)})$ characterize another class of solutions of
    the inverse problem ($c>0$ a norming constant) which is really a new one.
  \item Let again $(\mu_{k},m_{k})_{k\in\N}$ be a solution of the inverse
    problem in the type $II_{\infty}$ case, and let $k,l\in\N$ be a pair of
    indices and $\epsilon>0$. Then we get another class by adding $\epsilon$
    to $m_{k}$ and subtracting it from $m_{l}$ where again we have really a
    new class.
  \end{enumerate}
\end{bsp}

\begin{rem}\label{rem5:2}
\begin{enumerate}
  \item \bspref{bsp5:1}.3 and \bspref{bsp5:1}.4 shows that in the type
    $II_{\infty}$ case, when $\Op{H}_{0}$ has pure point spectrum, we can
    always construct a second class of solutions, different from the simple
    class discussed in \secref{sec:4},
    i.e. $NF_{\M_{0}}\not=NF_{\M_{0}}^{1}$, in contrast to
    the type $I$ case, where for modular operators with generic spectrum we
    have $NF_{\M_{0}}=NF_{\M_{0}}^{1}$ (cf. \cite{WolII}).
  \item
    Unfortunately the classification result presented here applies only to
    operators with pure point spectrum. Whereas in general there are also
    operators with more complicated spectrum (cf. \cite[Remark 4.1]{BolI}),
    for type $I$ factors this is no restriction, since all operators
    generating modular operators are trace class operators, hence have pure
    point spectrum. 
  \end{enumerate}
\end{rem}

\textbf{Acknowledgements}

The author thanks professor M. Wollenberg for discussing and his usefull hints
and the DFG and the Graduiertenkolleg for financial support.

\nocite{WolI}

\bibliographystyle{alpha}

\end{document}